\documentclass[12pt]{article}

\newcommand{\be}{\begin{equation}}
\newcommand{\ee}{\end{equation}}
\newcommand{\bi}[1]{\vspace{-3mm} \bibitem{#1}}

\begin{document}

{\it International Journal of Applied Mathematics.
Vol.27. No.5. (2014) 491-518. }
\vskip 3mm

\begin{center}
{\Large \bf FRACTIONAL-ORDER VARIATIONAL DERIVATIVE}
\vskip 5 mm

{\large \bf Vasily E. Tarasov }\\

\vskip 3mm

{\it Skobeltsyn Institute of Nuclear Physics, \\
Lomonosov Moscow State University, \\
Moscow 119991, RUSSIA} \\
\vskip 3 mm
e-mail: tarasov@theory.sinp.msu.ru
\end{center}

\begin{abstract}
We consider some possible approaches to
the fractional-order generalization of definition of 
variation (functional) derivative.
Some problems of formulation of a fractional-order
variational derivative are discussed.
To give a consistent definition of 
the fractional-order variations, 
we use a fractional generalization of Gateaux differential.
\end{abstract}

\noindent
%%%MSC 
Math. Subj. Classification 2010: 26A33; 49S05 
%%% 26A33 Fractional derivatives and integrals
%%% 49S05 Variational principles of physics 

%%%PACS: 45.10.Hj; 45.20.-d; 45.20.Jj; 
%%% 45.10.Hj Perturbation and fractional calculus methods
%%% 45.20.-d Formalisms in classical mechanics 
%%% 45.20.Jj Lagrangian and Hamiltonian mechanics 


\section{Introduction}

Derivatives of non-integer order are well-known in mathematics
\cite{SKM}-\cite{FC2}. 
The fractional calculus has a long history
from 1695, when the derivative of
order $\alpha=0.5$ has been described by Leibniz 
\cite{His1,His2,His3,His4,His5}.
Derivatives and integrals of fractional order have found many
applications in recent studies in mechanics and physics.
The interest to fractional equations has been growing continually 
during the last few years because of numerous applications 
\cite{Mainardi-Book}-\cite{APSS2014a}.

In mathematics and theoretical physics, 
the variational (functional) derivative is a 
generalization of the usual derivative 
that arises in the calculus of variations. 
In a variational (functional) derivative, instead of 
differentiating a function with respect to a variable, 
one differentiates a functional with respect to a function. 
The fractional generalization of variational (functional) 
derivative was suggested in this paper. 

In this paper, some problems of formulation of 
a fractional-order variational derivative are discussed.
To give a definition of fractional variation, 
we suggest to use a fractional generalization
of Gateaux differential.

In Section 2, some properties of the Riemann-Liouville or 
Caputo fractional derivatives are noted.
In Section 3, we give definitions of 
variational (functional) derivatives of integer orders.
In Section 4, we discuss problems of 
different possible ways to define a fractional generalization 
of variational (functional) derivatives.
In Section 5, we suggest a definition of 
fractional-order variational (functional) derivatives
by using the proposed fractional generalization
of Gateaux differential.
In Section 6, a fractional variation of fields 
that is defined by fractional exterior derivatives
are considered.
A conclusion is given in Section 7.

\section{Fractional Derivative}

The theory of derivatives of non-integer order goes back 
to Leibniz, Liouville, Riemann, Gr\"unwald, and Letnikov 
\cite{SKM}-\cite{FC2}. 
The authors of many papers use 
the fractional derivative $D^{\alpha}_{x}$
in the Riemann-Liouville or Caputo forms \cite{SKM}. 
Let us give definitions of these derivatives
and some properties. \\

{\bf Definition 1.} {\it 
The Riemann-Liouville fractional derivative 
of the function $f(x)$ belonging to 
the space $AC^n[a,b]$ of absolutely continuous functions
is defined on $[a,b]$ by the equation
\be \label{df} 
D^{\alpha}_{x}f(x)=\frac{1}{\Gamma (m-\alpha)}
\frac{d^m}{d x^m} 
\int^x_{a} \frac{f(y) \, dy}{(x-y)^{\alpha-m+1}} , \ee
where $\Gamma (z)$ is the Gamma function,
$m$ is the first whole number greater than or equal to $\alpha$.} \\

In equation (\ref{df}), the initial point of the fractional 
derivative \cite{SKM} is set to zero.
The derivative of powers $k$ of $x$ is
\be \label{xk} D^{\alpha}_x (x-a)^k=
\frac{\Gamma(k+1)}{\Gamma(k+1-\alpha)} x^{k-\alpha} ,
\quad (x>a) ,
\ee
where $k \ge 1$, and $\alpha \ge 0$. 
Note that the derivative of a constant $C$ need not be zero:
\[ D^{\alpha}_{x_i} C=\frac{x^{-\alpha}_i}{\Gamma(1-\alpha)} C . \]
Therefore we see that constants $C$ in the equation
$V(x)=C$ cannot define a stationary state for the equation
$D^{\alpha}_x V(x)=0$. 
In order to define stationary values, 
we should consider solutions of the equations 
$D^{\alpha}_{x_i} V(x)=0$.

The Riemann-Liouville fractional derivative has some notable 
disadvantages in physical applications such as 
the hyper-singular improper integral, where the order 
of singularity is higher than the dimension, and 
nonzero of the fractional derivative of constants, 
which would entail that dissipation does not vanish 
for a system in equilibrium. 
The desire to formulate initial value problems for physical
systems leads us to the use of Caputo fractional derivatives 
\cite{Caputo,Caputo2,Mainardi} 
(see also \cite{Podlubny,KST})
rather than Riemann-Liouville fractional derivative. \\

{\bf Definition 2.} {\it 
The Caputo fractional derivative 
of the function $f(x)$ belonging to 
the space $AC^n[a,b]$ of absolutely continuous functions
is defined on $[a,b]$ by the equation
\be \label{Caputo} 
D^{\alpha}_{x}f(x)=\frac{1}{\Gamma (m-\alpha)}
\int^x_{a} \frac{dy}{(x-y)^{\alpha-m+1}} 
\frac{d^m f(y)}{d y^m}, 
\quad (x>a) ,
\ee
where $f^{(m)}(y)=d^m f(y)/dy^m$, and
$m$ is the first whole number greater than or equal to $\alpha$.} \\

This definition is of course more restrictive that (\ref{df}),
in that requires the absolute integrability of the derivative
of order $m$. 
The Caputo fractional derivative first computes an ordinary
derivative followed by a fractional integral to achieve the
desire order of fractional derivative.
The Riemann-Liouville fractional derivative 
is computed in the reverse order. 
Integration by part of (\ref{Caputo}) leads us to 
\be \label{C-RL}
D^{\alpha}_{* \ x}f(x)=D^{\alpha}_{x}f(x)-
\sum^{m-1}_{k=0}
\frac{x^{k-\alpha}}{\Gamma(k-\alpha+1)} f^{(k)}(0+) .\ee 
It is observed that the second term in Eq. (\ref{C-RL}) regularizes 
the Caputo fractional derivative to avoid the potentially divergence 
from singular integration at $x=0+$. In addition, the 
Caputo fractional differentiation of a constant results in zero.
If the Caputo fractional derivative is used instead of the 
Riemann-Liouville fractional derivative then 
the stationary values 
are the same as those for the usual case ($V(x)-C=0$). 
The Caputo formulation of fractional calculus can be more 
applicable to definition of fractional variation 
than the Riemann-Liouville formulation.

%%%%%%%%%%%%%%%%%%%%%%%%%%%%%%%%%%%%
\section{Variational Derivatives of Integer Order}

In mathematics and theoretical physics, 
the variational (functional) derivative is a 
generalization of the usual derivative 
that arises in the calculus of variations. 
In a variational (functional) derivative, 
instead of differentiating a function 
with respect to a variable, one differentiates a functional 
with respect to a function.

%%%The possible definitions suitable for certain 
%%%computations is given here. 
%%%There are more general definitions of functional derivatives.

\subsection{Definition by increment and Taylor series}

%%%{\bf (1)} 
The variational derivative can be defined by the following way.
Let us consider the increment of the functional $F[u]$ that
is defined by the equation
\be
\Delta F[u]=F[u+h]-F[u] ,
\ee
and consider an integer-order variational derivative. \\
%%%=\int^{x_2}_{x_1} f(x,u+h) dx-
%%%\int^{x_2}_{x_1} f(x,u) dx . 
%%%\ee
%%%\be \Delta F[u]=F[u+\delta u]-F[u]=
%%%\int^{x_2}_{x_1} f(x,u+\delta u) dx-
%%%\int^{x_2}_{x_1} f(x,u) dx . \ee

{\bf Definition 3.} {\it 
If this increment of the functional $F[u]$ exists, and
can be represented in the form
\be
\Delta F[u]=\delta F(u,h)+\omega(h,u) ,
\ee
where 
\[\lim_{||h|| \rightarrow 0} \frac{||\omega(h,u)||}{||h||}=0 ,\]
then $\delta F$ is called the first variation or 
Frechet derivative \cite{F1,F2,Vainberg} 
of functional $F$.
The function $h=h(x)$ is called the variation, 
and it is denoted by $\delta u$.} \\

%%%%%%%%%%%%%%%%%%%%%%%%%%%%%%%%%%%%%%%%%%%%%%%%%%%%

{\bf Example 1.} Let us define the functional
\be \label{Fun0}
F[u]=\int^{x_2}_{x_1} f(u) dx
\ee
in some Bahan space $E$.
The increment of the functional $F[u]$ 
is defined by the equation
\be
\Delta F[u]=F[u+h]-F[u]=
\int^{x_2}_{x_1} f(u+h) dx-
\int^{x_2}_{x_1} f(u) dx . 
\ee
%%%\be \Delta F[u]=F[u+\delta u]-F[u]=
%%%\int^{x_2}_{x_1} f(x,u+\delta u) dx-
%%%\int^{x_2}_{x_1} f(x,u) dx . \ee
Here we suppose that the variation $h(x)=\delta u(x)$
is equal to zero in boundary points $x_1$ and $x_2$. 
Expanding the integrand $f(x,u+h)$
in the power series up to first order
%%% fractional Taylor
\[ f(u+h)=f(u)+ \frac{\partial f(u)}{\partial u} h+R(h,u). \]
where 
\be
\lim_{h\rightarrow 0} h^{-1} R(h,u)=0 .
\ee
%%% Linear Property.
Therefore we have 
\be
\delta F[u]=
%%%F[u+\delta u]-F[u]=
\int^{x_2}_{x_1} \frac{\partial f(u)}{\partial u} h \, dx . 
\ee
The variational derivative for functional (\ref{Fun0})
has the form
\be \frac{\delta F[u]}{\delta u}=\frac{\partial f(u)}{\partial u} . \ee

%%%%%%%%%%%%%%%%%%%%%%%%%%%%%%%%%%%%%%%%%%%%%%%%%%%%

{\bf Example 2.} 
Let us consider the functional
\be \label{Fun1}
F[u]=\int^{x_2}_{x_1} f(u,u_x) dx .
\ee
We can derive the first variation of the functional by the equation
\be
\delta F[u]=F[u+\delta u]-F[u]=
\int^{x_2}_{x_1} f(u+h,[u+h]_x) dx-
\int^{x_2}_{x_1} f(u,u_x) dx . 
\ee
Here we suppose that the variation $h(x)=\delta u(x)$
is equal to zero in boundary points $x_1$ and $x_2$.
Let us expand the integrand $f(u+\delta u,(u+\delta u)_x)$
in the power series up to first order
%%% fractional Taylor
\[ f(u+h,(u+h)_x)=f(u,u_x)+
\frac{\partial f}{\partial u} h+
\frac{\partial f}{\partial u_x} h_x , \]
where $h_x=dh/dx$.
Then we get the variation of the functional
\be
\delta F[u]=\int^{x_2}_{x_1} \left[ 
\frac{\partial f}{\partial u} h+
\frac{\partial f}{\partial u_x} h_x \right] dx . \ee
Integrating the second term by part, and supposing 
\[ h(x_1)=h(x_2)=0 , \]
we get the result
\be
\delta F[u]=\int^{x_2}_{x_1} \left[ \frac{\partial f}{\partial u} -
\frac{d}{dx} 
\left(\frac{\partial f}{\partial u_x} \right) \right] h(x) \, dx . \ee
Using this relation, we get
the variational derivative of the functional (\ref{Fun1})
in the form
\be
\frac{\delta F[u]}{\delta u}=\frac{\partial f}{\partial u} -
\frac{d}{dx} \left(\frac{\partial f}{\partial u_x} \right) . 
\ee

%%%%%%%%%%%%%%%%%%%%%%%%%%%%%%%%%%%%%%%%%%%%%%%%%%%%

{\bf Example 3.} If we consider the functional
\be
F[u]=\int^{x_2}_{x_1} f(x,u,u_x,...u^{(n)}) dx
\ee
then the variation of the functional
is defined by the equation
\be
\delta F[u]=\int^{x_2}_{x_1} \sum^{n}_{m=0} 
\frac{\partial f}{\partial u^{(m)}} (\delta u)^{(m)} dx ,
\ee
and the variational derivative has the form
\be
\frac{\delta F[u]}{\delta u}=\sum^{n}_{m=0} (-1)^m
\frac{d^m}{dx^m} \left(\frac{\partial f}{\partial u^{(m)}} \right) , 
\ee
where $u^{(m)}=d^m u(x)/dx^m$.

\subsection{Definition by composite function of parameter}

There is the other approach to definition of variation 
and variational derivative.
Let us consider the functional (\ref{Fun1}),
where $u$ is a function of coordinates $x$ and parameter $a$,
i.e., $u=u(x,a)$, and
\be
F[u]=\int^{x_2}_{x_1} f(u(x,a),u_x(x,a)) dx .
\ee
The derivative of $F$ with respect $a$ can be derived
%%% differential !!!
\be
\frac{d F}{da}=\int^{x_2}_{x_1}dx \frac{df}{da} =\int^{x_2}_{x_1}\left[ 
\frac{\partial f}{\partial u} \frac{\partial u}{\partial a}+
\frac{\partial f}{\partial u_x} \frac{\partial u_x}{\partial a}
\right] dx .
\ee
Using
\[ \frac{\partial u_x}{\partial a}=\frac{\partial}{\partial a} 
\frac{\partial u}{\partial x}=\frac{\partial}{\partial x} 
\frac{\partial u}{\partial a}, \]
the conditions $\delta u(x_1,a)=\delta u(x_2,a)=0$,
and integrating by part, we get
\be
\frac{d F}{da}=\int^{x_2}_{x_1} 
\frac{\partial u}{\partial a} 
\left[ \frac{\partial f}{\partial u} -
\frac{\partial}{\partial x} \left(\frac{\partial f}{\partial u_x} \right)
\right] dx . \ee
As the result, we have
\be
\frac{d F}{da}=\int^{x_2}_{x_1} 
\frac{\partial u(x)}{\partial a} 
\frac{\delta F[u]}{\delta u(x)} dx . \ee

A fractional-order generalization of this approach 
is difficult to realize. 
The reason of this difficulty is connected with
difficulty to generalized the rule of 
differentiating a composite functions
\be
\frac{df(u(a))}{da}=\frac{\partial f}{\partial u} 
\frac{\partial u}{\partial a} 
\ee
for the fractional case.
It is known that coordinate transformations is connected with the derivative of a composite function 
$D^1_t f(u(x))= (D^{1}_uf)(u=u(x)) \, (D^1_x u)(x)$.
The formula of fractional derivative of a composite function 
(see Equation 2.209, Section 2.7.3 page 98 \cite{Podlubny})
is the following
\be
D^{\alpha}_x f(u(x)) = \frac{(x-a)^{\alpha}}{\Gamma(1-\alpha)} 
f(u(x))
+\sum^{\infty}_{k=1} C^{\alpha}_k \, \frac{k! (x-a)^{k-\alpha}}{\Gamma(k-\alpha+1)} \,
\sum^k_{m=1} (D^m f)(u(x)) \sum \prod^k_{r=1} \frac{1}{a_r!} \Bigl(\frac{(D^r_xu)(x))}{r!}\Bigr)^{a_r}
\ee
where the sum $\sum$ extends over all combinations of 
non-negative integer values of $a_1$, $a_2$, . . . , $a_k$ 
such that 
\be
\sum^k_{r=1} r a_r=k , \quad \sum^k_{r} a_r=m.
\ee

%%%%%%%%%%%%%%%%%%%%%%%%%%%%%%%%%%%%%%%%%%%%%%%%%%%%%%%%%%%%%%
\section{Problems to Formulate 
Fractional-Order Variational Derivative}

Let us consider some problems in the formulation
of the fractional-order generalization
of variational derivatives.

%%%In calculus, a differential is an infinitesimal change 
%%%in the value of a function. 
%%%Here we suppose that the variation $\delta u(x)$
%%%is equal to zero in boundary points $x_1$ and $x_2$.

{\bf Problem 1.}
To define the fractional variation, we can use
the fractional Taylor series 
(see Section 2.6 in \cite{SKM}).
Expanding the integrand $f(u+h)$
in the fractional Taylor power series 
\be f(u+h)=
\frac{1}{\Gamma(\alpha+1)} (h)^{\alpha} D^{\alpha}_u f(u) +
\frac{1}{\Gamma(\alpha+2)} (h)^{\alpha+1} D^{\alpha}_u f(u)+
... ,\ee
we can see that this series cannot have the term $f(u)$.
Therefore, we cannot consider the increment
\be
\Delta f(u)=f(u+h)- f(u).
\ee
In order to use Taylor series for the definition, 
we can consider the increment in the form
\be
\Delta_q f=f(u+h)- f(u+qh).
\ee
where $0<q<1$. In this case, we have
\be 
\Delta f(u)=
\frac{1-q^{\alpha}}{\Gamma(\alpha+1)} h^{\alpha} D^{\alpha}_u f(u) +
\frac{1-q^{\alpha+1}}{\Gamma(\alpha+2)} h^{\alpha+1} D^{\alpha}_u f(u)
+ ... .\ee
This approach can be connected with the q-analysis \cite{Kac},
and fractional q-derivatives \cite{FQ1,FQ2}.
We will consider the fractional variational q-derivative 
in the next paper. 

We can consider the increment of the functional $F[u]$ that
is defined by the equation
\be
\Delta_q F[u]=F[u+h]-F[u+qh] . \ee
%%%=\int^{x_2}_{x_1} f(x,u+h) dx-
%%%\int^{x_2}_{x_1} f(x,u) dx . 
%%%\ee
%%%\be \Delta F[u]=F[u+\delta u]-F[u]=
%%%\int^{x_2}_{x_1} f(x,u+\delta u) dx-
%%%\int^{x_2}_{x_1} f(x,u) dx . \ee
If this increment of the functional $F[u]$ exists, and
can be represented in the form
\be
\Delta_q F[u]=\delta^{\alpha} F(u,h^{\alpha})+\omega(h^{\alpha},u) ,
\ee
where 
\[\lim_{||h|| \rightarrow 0} 
\frac{||\omega(h^{\alpha},u)||}{||h^{\alpha}||}=0 ,\]
then $\delta^{\alpha} F$ can be considered as a fractional variation 
of functional $F$. \\

%%%%%%%%%%%%%%%%%%%%%%%%%%%%%%%%%%%%%%%%%%%%%%%%%%%%%%%%%%%%%

{\bf Problem 2.}
 To define the fractional generalization of
variation and fractional exterior variational calculus \cite{AK},
we can use an analogy with definition of fractional 
exterior derivative.
If the partial derivatives in the definition of the exterior derivative 
\[ d=dx_i \partial / \partial x_i \] 
are allowed to assume fractional order, 
a fractional exterior derivative can be defined \cite{FDF} 
by the equation
\be
d^{\alpha}=(dx_i)^{\alpha} D^{\alpha}_{x_i} , 
\ee 
where $D^{\alpha}_{x}$ are the fractional derivative with respect to $x$.
Using this analogy, we can define the fractional variation
by the following way.
For the point $u$ of functional space, 
we can define the fractional variation
$\delta F[u]$ of the functional
\[ F[u]=\int^{x_2}_{x_1} f(u,u_x) dx , \]
where $u_x=du/dx$, by the equation
\be
\delta^{\alpha} F[u]= \int^{x_2}_{x_2} dx \left[ (\delta u)^{\alpha}
D^{\alpha}_u f(u,u_x) +
(\delta u_x)^{\alpha} D^{\alpha}_{u_x} f(u,u_x) \right].
\ee 
Unfortunately, this approach leads to the difficulties
with the realization of integration by part in the second term.
It is easy to see that the variation $(\delta u_x)^{\alpha}$
cannot be represented as some operator $\hat{A}(d/dx)$ acts on the variation
$(\delta u)^{\alpha}$, i.e., we have
\be \label{dif}
\left( \delta u_x(x) \right)^{\alpha}= 
\left(\frac{d}{dx} \delta u(x) \right)^{\alpha} 
\not= \hat{A}\left(\frac{d}{dx}\right) (\delta u)^{\alpha}.
\ee
In the particular case,
$\left( \delta u_x(x) \right)^{\alpha}= 
\left(d(\delta u)/dx \right)^{\alpha} 
\not= d (\delta u)^{\alpha}/dx$ if $\alpha \ne 1$.

%%%%%%%%%%%%%%%%%%%%%%%%%%%%%%%%%%%%%%%%%%%%%%%%%%%%%%%%%%%%%%%%

{\bf Problem 3.}
To define the fractional variational derivative, 
we can use the fractional derivative with respect to function 
(see Section 18.2 in \cite{SKM}).
The fractional Riemann-Liouville defivative of the function $f(x)$ 
with respect the function $u(x)$ of order $\alpha$ ($0<\alpha<1$):
\be
D^{\alpha}_{0+,u(x)} f(x)=
\frac{1}{\Gamma(1-\alpha)} \left( \frac{d u(x)}{dx}\right)^{-1}
\frac{d}{dx} \int^x_0 dy \frac{f(y)}{[u(x)-u(y)]^{\alpha}} \frac{d u(y)}{dy}.
\ee
Using the function $u=u(x,a)$, we can define
fractional derivative with respect to function of parameter
\be
D^{\alpha}_{0+,u(x,a)} f(x,a)=
\frac{1}{\Gamma(1-\alpha)} \left( \frac{d u(x,a)}{da}\right)^{-1}
\frac{d}{da} \int^a_0 d b \frac{f(x,b)}{[u(x,a)-u(x,b)]^{\alpha}} 
\frac{d u(x,b)}{db}. \ee
As the result, the fractional variation can be defined by the equation
\be
\delta^{\alpha} F[u(x)]=
\int^{x_2}_{x_1} dx D^{\alpha}_{u(x,a)} f(u(x,a)) (\delta u(x,q))^{\alpha} ,
\quad \frac{\delta^{\alpha} F[u]}{\delta u^{\alpha}}= 
D^{\alpha}_{u(x,a)} f(u(x)) .
\ee
This approach leads to the difficulties with the realization 
of integration by part in the second term (\ref{dif}).
In this paper, this approach to definition is not considered.

To avoid all these difficulties and problems,
we suggest to define a fractional variational derivative
by using a fractional-order generalization of Gateaux differential.

\section{Fractional Generalization of Gateaux Differential}

In this section, we consider 3 steps to define the fractional 
generalization of variation by using some generalization
of Gateaux differential.
We consider these steps in order to explain the final definition
of fractional Gateaux variation.

\subsection{Variations of integer order}

Suppose the functional $F[u]$ 
is continuous (smooth) map (with certain boundary conditions) 
from everywhere dense subset $D(F)$ of Banah space to space $R$. 
Let us define the Gateaux differential \cite{G,F1,F2} of 
a functional $F[u]$ at the point $u(x)$ of subset $D(F)$ 
of the functional Banah space. \\

{\bf Definition 4.} {\it 
The Gateaux differential (or first variation)
of a fucntional $F[u]$ is defined by the equation
\be \label{VD0} 
\delta F[u,h] =
\left( \frac{d}{d \epsilon} F[u+\epsilon h] \right)_{\epsilon=0}=
\lim_{\epsilon \rightarrow 0} 
\frac{F[u +\epsilon \ h]-F[u]}{\epsilon} \ee
if the limit exists for all $h(x)\in D(F)$. 
The function $h(x)$ is called a variation of function $u(x)$ 
and denoted by $\delta u(x)=h(x)$.} \\

A first variation of functional $F[u]$ at the point $u=u(x)$ 
is defined as a first derivative of functional $F[u+\epsilon h]$
with respect to parameter $\epsilon$ for $\epsilon=0$.

If the Gateaux's differential of the functional
\be F[u]=\int dx f(x,u,u_x,...) \ee
is a linear with respect to $h(x)$, then we can write
\be
\delta F[u,h]=\int dx \, E(x,u,u_x,...) h(x) , 
\ee
and $E(x,u,u_x,...)$ is called the variational derivative
and is denoted by $ \delta F/ \delta u$. \\

%%%For a functional $F$ mapping 
%%%(continuous (smooth) map with certain boundary conditions)
%%%functions $u$ from a manifold $M$ to $R$ or $C$, 
%%%the variation (functional) derivative of $F$
%%%(the Gateaux differential \cite{G,F1,F2} of a functional $F$) 
%%%, denoted $\delta F$ 
%%%is a distribution such that for all test functions $f$,
%%%\be \label{VD0} \delta F[u] =
%%%\left( \frac{d}{d \epsilon} F[u+\epsilon h] \right)_{\epsilon=0}=
%%%\lim_{\epsilon \rightarrow 0} 
%%%\frac{F[u +\epsilon \ h]-F[u]}{\epsilon}. \ee
%%%
%%%A first variation of functional $F[u]$ at the point $u=u(x)$ 
%%%can be defined as a first derivative of functional $F[u+\epsilon h]$
%%%with respect to parameted $\epsilon$ for $\epsilon=0$:
%%%\be \label{VD0}
%%%\delta F[u]=\left( \frac{d}{d\epsilon} 
%%%F[u+\epsilon h]\right)_{\epsilon=0}. \ee
%%%The variation of function $u(x)$ is defined 
%%%by $\delta u(x)=h(x)$.

{\bf Example.}
For the example, we can consider the functional
\be \label{u2}
F[u]=\int^{x_2}_{x_1} [u(x)]^2 dx .
\ee
The functional $F[u+\epsilon h]$ has the form
\be \label{F}
F[u+\epsilon h]=\int^{x_2}_{x_1} [u(x)]^2 dx+
2\epsilon \int^{x_2}_{x_1} u(x) h(x) dx+
\epsilon^2\int^{x_2}_{x_1} [h(x)]^2 dx .
\ee
The derivative of first order with respect to parameter
$\epsilon$ is
\be
\frac{d}{d\epsilon} F[u+\epsilon h] =2\int^{x_2}_{x_1} u(x) h(x) dx+
2 \epsilon \int^{x_2}_{x_1} [h(x)]^2 dx .
\ee
As the result, the variational derivative is equal
\be
\left(\frac{d}{d\epsilon} F[u+\epsilon h]\right)_{\epsilon=0}
=2\int^{x_2}_{x_1} u(x) h(x) dx .
\ee

\subsection{STEP 1.}

It seems that we can define a fractional variation by the equation
\be \label{fVD0} 
\delta^{\alpha} F[u,h] =
\left( \left[\frac{d}{d \epsilon}\right]^{\alpha} 
F[u+\epsilon h] \right)_{\epsilon=0}. \ee

Let us consider the fractional Riemann-Liouville derivative of
functional (\ref{F}) with respect to parameter $\epsilon$:
\be 
D^{\alpha}_{\epsilon}F[u+\epsilon h]=
(D^{\alpha}_{\epsilon}1)\int^{x_2}_{x_1} [u(x)]^2 dx+
2(D^{\alpha}_{\epsilon}\epsilon )\int^{x_2}_{x_1} u(x) h(x) dx+
(D^{\alpha}_{\epsilon}\epsilon^2)\int^{x_2}_{x_1} [h(x)]^2 dx .
\ee
Using equation
\be
D^{\alpha}_{\epsilon} \epsilon^k=\frac{\Gamma(k+1)}{\Gamma(k+1-\alpha)} 
\epsilon^{k-\alpha},
\ee
we have
\[ D^{\alpha}_{\epsilon}F[u+\epsilon h]=
\frac{\epsilon^{-\alpha}}{\Gamma(1-\alpha)} 
\int^{x_2}_{x_1} [u(x)]^2 dx+
\frac{2 \epsilon^{1-\alpha}}{\Gamma(2-\alpha)} 
\int^{x_2}_{x_1} u(x) h(x) dx+ \]
\be \label{DF}
+\frac{2}{\Gamma(3-\alpha)} 
\epsilon^{k-\alpha}\int^{x_2}_{x_1} [h(x)]^2 dx .
\ee

If the derivatives with respect to parameter $\epsilon$ 
in the definition of variational (\ref{VD0})
are allowed to assume fractional order, 
a fractional variational can be defined.
Unfortunately, if we define the fractional variation of 
the functional by the equation
\be \label{FV1}
\delta^{\alpha} F[u]=\left( D^{\alpha}_{\epsilon} 
F[u+\epsilon h]\right)_{\epsilon=0},
\ee
then we have some problems that state about incorrectness
of definition (\ref{FV1}). These problems are following:

1) The first term of right hand side of equation (\ref{DF}) 
leads us to infinity. If $\epsilon$ tends to zero, 
then we get $\epsilon^{-\alpha}\rightarrow \infty$,
and $D^{\alpha}_{\epsilon}F[u+\epsilon h] \rightarrow \infty$.
Therefore the first term that is follows from the relation 
$D^{\alpha}_{\epsilon} C=0$ must be removed. 
For this aim we can use the Caputo fractional derivative 
\cite{Podlubny,Caputo,Caputo2,Mainardi}.
Using (\ref{C-RL}) for $0<\alpha<1$, we have
\be \label{C-RL-1}
D^{\alpha}_{* \ \epsilon} f(\epsilon)=
D^{\alpha}_{\epsilon}f(\epsilon)-
\frac{\epsilon^{-\alpha}}{\Gamma(2-\alpha)} f(0+) .\ee 

2) The second term of the right hand side of equation (\ref{DF})
is proportional to $\epsilon^{1-\alpha}$.
This proportionality leads us to zero in the limit $\epsilon \rightarrow 0$,
and we cannot derive some nonzero relation. Therefore we must consider
the functional $F[u+(\epsilon h)^{\alpha}]$ in the definition. \\

As the result, we cannot use the definition (\ref{fVD0}) and (\ref{FV1}).

\subsection{STEP 2.}

Let us consider the following equations 
for the definition of fractional variational functional $F[u]$:
\be \label{FV2}
\delta^{\alpha} F[u]=\left( D^{\alpha}_{\epsilon} 
F[u+(\epsilon h)^{\alpha}] \right)_{\epsilon=0},
\ee
where $D^{\alpha}_{\epsilon}=D^{\alpha}_{*}$ is a Caputo fractional 
derivative of order $\alpha$ with respect to $\epsilon$.

If we consider the functional (\ref{u2}), then
\be 
F[u+(\epsilon h)^{\alpha}]=\int^{x_2}_{x_1} [u(x)]^2 dx+
2\epsilon^{\alpha} \int^{x_2}_{x_1} u(x) h^{\alpha}(x) dx+
\epsilon^{2\alpha}\int^{x_2}_{x_1} [h(x)]^{2\alpha} dx .
\ee
The fractional derivative (\ref{FV2}) of this functional is equal to
\be \label{DF2}
D^{\alpha}_{\epsilon}F[u+(\epsilon h)^{\alpha}]=
2 \Gamma(\alpha+1) \int^{x_2}_{x_1} u(x) h^{\alpha}(x) dx+
\frac{\Gamma(2\alpha+1)}{\Gamma(\alpha+1)} 
\epsilon^{\alpha}\int^{x_2}_{x_1} [h(x)]^{2\alpha} dx .
\ee
Therefore, we get the following equation
\be \label{DF3}
\left(D^{\alpha}_{\epsilon}F[u+(\epsilon h)^{\alpha}]\right)_{\epsilon=0}=
2 \Gamma(\alpha+1) \int^{x_2}_{x_1} u(x) h^{\alpha}(x) dx ,
\ee
%%%or, in an equivalent form
%%%\be
%%%\delta^{\alpha} F[u]=
%%%2 \Gamma(\alpha+1) \int^{x_2}_{x_1} u(x) ({\delta u}(x) )^{\alpha} dx .
%%%\ee
%%%The variational derivative of fractional order $\alpha$ can be defined
%%%\be
%%%\frac{\delta^{\alpha} F[u]}{{\delta u}^{\alpha}(x)}=
%%%2 \Gamma(\alpha+1) u(x) . \ee
Note that for $\alpha=1$, we get the usual relation
\be
\delta^{\alpha=1} F[u]=2 \int^{x_2}_{x_1} u(x) h(x) dx .
\ee
Unfortunately, if we use $\alpha=0$, we get 
\be
F[u]=\delta^{\alpha=0} F[u]=
2 \int^{x_2}_{x_1} u(x) (h (x))^{0} dx =
2 \int^{x_2}_{x_1} u(x) dx \not=F[u].
\ee
Therefore we cannot use the definition (\ref{FV2}).

\subsection{STEP 3.}

The usual definition (\ref{VD0}) can be rewritten in the form
\be \label{Var}
\delta F[u]=\left( \frac{d}{d\epsilon} 
F\left[u \left(1+\frac{\epsilon h}{u} \right)\right]\right)_{\epsilon=0}.
\ee
The fractional variation can be defined by the equation
\be \label{FV3}
\delta^{\alpha} F[u]=\left( D^{\alpha}_{\epsilon}
F\left[u\left(1+\left[\frac{\epsilon h}{u}\right]^{\alpha}\right)\right]
\right)_{\epsilon=0}.
\ee
This definition is more consistent in order to
realize the physical dimensions. In this definition
the parameter $\epsilon$ is dimensionless.

If we consider the functional (\ref{u2}), then
\[ F\left[u\left(1+\left[\frac{\epsilon h}{u}\right]^{\alpha}\right)\right]=
F[u + (\epsilon h)^{\alpha} u^{1-\alpha}]= \]
\be 
=\int^{x_2}_{x_1} [u(x)]^2 dx+
2\epsilon^{\alpha} \int^{x_2}_{x_1} u^{2-\alpha}(x) h^{\alpha}(x) dx+
\epsilon^{2\alpha}\int^{x_2}_{x_1} u^{2-2\alpha}(x) [h(x)]^{2\alpha} dx .
\ee
The fractional derivative of this functional is equal to
\be 
D^{\alpha}_{\epsilon} F[u + (\epsilon h)^{\alpha} u^{1-\alpha}]=
2 \Gamma(\alpha+1) \int^{x_2}_{x_1} u^{2-\alpha}(x) h^{\alpha}(x) dx+
\frac{\Gamma(2\alpha+1)}{\Gamma(\alpha+1)} 
\epsilon^{\alpha}\int^{x_2}_{x_1} u^{2-2\alpha}(x) [h(x)]^{2\alpha} dx .
\ee
As the result, we get
\be 
\left( D^{\alpha}_{\epsilon} F[u + (\epsilon h)^{\alpha} u^{1-\alpha}] \right)_{\epsilon=0}
=2 \Gamma(\alpha+1) \int^{x_2}_{x_1} u^{2-\alpha}(x) h^{\alpha}(x) dx .
\ee
For $\alpha=1$, we get the usual relation. 
For $\alpha=0$, we have 
\be 
F[u]=\delta^{\alpha=0} F[u]
=2 \int^{x_2}_{x_1} u^2(x) dx \not=F[u].
\ee
If we consider the functional
\be
F[u]= \int^{x_2}_{x_1} u^n(x) dx ,
\ee
then definition (\ref{FV3}) leads us to the relation
\be 
\left( D^{\alpha}_{\epsilon} 
F[u + (\epsilon h)^{\alpha} u^{1-\alpha}] \right)_{\epsilon=0}
=n \Gamma(\alpha+1) \int^{x_2}_{x_1} u^{n-\alpha}(x) h^{\alpha}(x) dx 
\ee
and, for $\alpha=0$
\be 
F[u]=\delta^{\alpha=0} F[u]
=n \int^{x_2}_{x_1} u^n(x) dx \not=F[u].
\ee
Therefore we must realize some changes in 
the definition (\ref{FV3}).
This modification is suggested in the next subsection.

%%%{\bf STEP 4.}

%%%%%%%%%%%%%%%%%%%%%%%%%%%%%%%%%%%%%%%%%%%%%%%%%%%%%%%%%%

\subsection{Definition of variation of fractional order}

Taking into account the remarks in the form of Step 1-3 and
using Eq. (\ref{Var}) for variation of integer order,
we can define the fractional variation in the following form. \\

{\bf Definition 5.} {\it 
The fractional-order variation $\delta^{\alpha} F[u]$
of the functional $F[u]$ is defined by the equation
\be \label{FV4}
\delta^{\alpha} F[u] =\left( D^{\alpha}_{\epsilon}
F\left[u\left(1+\frac{\epsilon h}{u}\right)^{\alpha} \right]
\right)_{\epsilon=0},
\ee
or, in an equivalent form
\be \label{FV5}
\delta^{\alpha} F[u]=\left( D^{\alpha}_{\epsilon}
F\left[ u^{1-\alpha} (u+ \epsilon h )^{\alpha} \right]
\right)_{\epsilon=0} ,
\ee
where $\epsilon \ge 0$.} \\

{\bf Example.} 
If we consider the functional (\ref{u2}), then we have
\be \label{exam}
F\left[ u^{1-\alpha} (u+ \epsilon h )^{\alpha} \right]= 
\int^{x_2}_{x_1} u^{2-2\alpha} h^{2\alpha} (\epsilon +u/h )^{2\alpha} dx.
\ee
Using the relation
\be
D^{\alpha}_{\epsilon} (\epsilon+u/h)^{2\alpha}=
\frac{\Gamma(2\alpha+1)}{\Gamma(\alpha+1)}(\epsilon+u/h)^{\alpha} ,
\ee
where we consider $\epsilon + u(x)/h(x) =\epsilon - \epsilon_0$
with $\epsilon_0=\epsilon_0(x)=-u(x)/h(x)$,
we get the fractional derivative of functional (\ref{exam}) in the form
\be 
D^{\alpha}_{\epsilon} F\left[ u^{1-\alpha} (u+ \epsilon h )^{\alpha} \right]= 
\frac{\Gamma(2\alpha+1)}{\Gamma(\alpha+1)} \int^{x_2}_{x_1} 
u^{2-2\alpha} h^{2\alpha} (\epsilon +u/h )^{\alpha} dx.
\ee
As the result, we get
\be 
\left( 
D^{\alpha}_{\epsilon} F\left[ u^{1-\alpha} (u+ \epsilon h )^{\alpha} \right]
\right)_{\epsilon=0}=
\frac{\Gamma(2\alpha+1)}{\Gamma(\alpha+1)} \int^{x_2}_{x_1} dx
u^{2-\alpha}(x) h^{\alpha}(x) .
\ee
For $\alpha=1$, and $\alpha=0$, we get the usual relations. 
Therefore definition (\ref{FV4}) satisfies the correspondent 
requirement for $\alpha=0$ and $\alpha=1$
As the result, we get the fractional variational derivative
of order $0\le \alpha \le 1$ in the form (\ref{FV5}). \\

{\bf Proposition 1.} {\it 
The fractional-order variation derivative (\ref{FV4})
of the functional 
\be \label{un}
F[u]= \int^{x_2}_{x_1} u^n(u) dx.
\ee
has the form
\be 
\delta^{\alpha} F [u]= \lambda (\alpha , n)
\int^{x_2}_{x_1} dx \, 
(D^{\alpha}_{u} u^n)(x) h^{\alpha}(x) ,
\ee
where 
\be
\lambda (\alpha , n) = 
\frac{\Gamma(n\alpha+1) \Gamma (n+1-\alpha)}{
\Gamma((n-1)\alpha+1) \Gamma (n+1)} .
\ee
} \\

{\bf Proof:} 
Using (\ref{un}), then we have
\be 
F\left[ u^{1-\alpha} (u+ \epsilon h )^{\alpha} \right]= 
\int^{x_2}_{x_1} u^{n-n\alpha} h^{n\alpha} (\epsilon +u/h )^{n\alpha} dx.
\ee
Using the relation
\be
D^{\alpha}_{\epsilon} (\epsilon+u/h)^{n\alpha}=
\frac{\Gamma(n\alpha+1)}{\Gamma((n-1)\alpha+1)}(\epsilon+u/h)^{(n-1)\alpha} ,
\ee
we get the fractional derivative of functional (\ref{un}) in the form
\be 
\delta^{\alpha} F [u]=
\left( 
D^{\alpha}_{\epsilon} F\left[ u^{1-\alpha} 
(u+ \epsilon h )^{\alpha} \right]
\right)_{\epsilon=0}=
\frac{\Gamma(n\alpha+1)}{\Gamma((n-1)\alpha+1)} 
\int^{x_2}_{x_1} dx u^{n-\alpha}(x) h^{\alpha}(x) .
\ee
This equation can be represented in the form 
\be \label{daDF}
\delta^{\alpha} F =
\int^{x_2}_{x_1} dx (D^{\alpha}_{u} u^n) h^{\alpha}(x) ,
\ee
up to numerical factor. For equation (\ref{daDF}), we have
\be
(D^{\alpha}_{u} u^n)= \frac{\Gamma(n+1)}{\Gamma(n+1-\alpha)} u^{n-\alpha} .
\ee
This end of proof. \\

%%%{\bf Definition 6.} {\it 
The fractional variational derivative
of order $\alpha \ge 1$ can be defined by the usual relation
\be
\delta^{\alpha} = \delta^{[\alpha]} \delta ^{\{\alpha\}}
\ee
by analogy with fractional derivative 
\be
D^{\alpha}_x=\frac{d^{[\alpha]}}{dx^{[\alpha]}} D^{\{\alpha\}}_x ,
\ee
%%%where $m$ is integer part of real number 
where $[\alpha]$ is a whole part of $\alpha$, and
$\{\alpha\}$ is a fractional part of number $\alpha$, i.e.,
$\{\alpha\}=\alpha-[\alpha]$.
%%%} 

%%%%%%%%%%%%%%%%%%%%%%%%%%%%%%%%%%%%%%%%%%%%%%%%%%%%%%%%%%%%%%%

\subsection{Functional with derivative of field}

Let us consider the functional $F[u]$ that has the form
\be \label{uux}
F[u]=\int^{x_2}_{x_1} f(u,u_x) dx
\ee
where $u_x=\frac{d u(x)}{dx}$. 
The first variation of this functional is defined by
\be \label{VD00}
\delta F[u]=\left( \frac{d}{d\epsilon} 
F[u+\epsilon h]\right)_{\epsilon=0} = 
\lim_{\epsilon \rightarrow 0}
\left( \frac{d}{d\epsilon} \int^{x_2}_{x_1} 
f(u+\epsilon h,[u+\epsilon h]_x) \right) dx ,
\ee
where we use
\be \label{ueh}
[u+\epsilon h]_x=u_x+\epsilon h_x .
\ee

In order to use the equation
\be
D^{\alpha}_{\epsilon} (\epsilon-\epsilon_0)^{\beta}=
\frac{\Gamma(\beta+1)}{\Gamma(\beta+1-\alpha)}
(\epsilon-\epsilon_0)^{\beta-\alpha} ,
\ee
for the relation $[\epsilon + u(x)/h(x)]$,
we consider $\epsilon_0=\epsilon_0(x)=-u(x)/h(x)$.
The definition of the fractional variation 
of functional (\ref{uux}) can be realized 
in the following form. \\

{\bf Definition 6.} {\it 
The fractional-order variation derivative (\ref{Def-6})
of the functional (\ref{uux}) is defined by the equation
\be \label{Def-6}
\delta^{\alpha} F[u] = \lim_{\epsilon \rightarrow 0}
\int^{x_2}_{x_1} \left( D^{\alpha}_{\epsilon}
f\Bigl( u^{1-\alpha} (u+ \epsilon h)^{\alpha},
[u^{1-\alpha} (u+ \epsilon h)^{\alpha} ]_x \Bigr)
\right) .
\ee
} \\

Using the definition (\ref{FV5}), we have 
the fractional variation of functional (\ref{uux}) in the form
\be \delta^{\alpha} F[u] =
%%%\left( D^{\alpha}_{\epsilon}
%%%F\left[u\left(1+\frac{\epsilon h}{u}\right)^{\alpha} \right]
%%%\right)_{\epsilon=0}=
\int^{x_2}_{x_1}\left( D^{\alpha}_{\epsilon}
f\Bigl( u^{1-\alpha} (u+ \epsilon h)^{\alpha},
[u^{1-\alpha} (u+ \epsilon h)^{\alpha} ]_x \Bigr)
\right)_{\epsilon=0} .
\ee
It is easy to see that the expression 
$[u^{1-\alpha} (u+ \epsilon h)^{\alpha} ]_x$
cannot be represented in the form of similar Eq. (\ref{ueh}).
This expression has the form
\be \label{ueh-f}
[u^{1-\alpha} (u+ \epsilon h)^{\alpha} ]_x=
[(1-\alpha)u^{-\alpha}(u+\epsilon h)^{\alpha}
+\alpha u^{1-\alpha} (u+\epsilon h)^{\alpha-1}] u_x+
\alpha u^{1-\alpha} (u+\epsilon h)^{\alpha-1} \epsilon h_x .
\ee
For $\alpha=0$, we have 
\be [u^{1-\alpha} (u+ \epsilon h)^{\alpha} ]_x=u_x, \ee
and for $\alpha=1$, we have
\be [u^{1-\alpha} (u+ \epsilon h)^{\alpha} ]_x=u_x+ \epsilon h_x .
\ee

Let us give the proposition for a special form of 
the functional. \\

%%%%%%%%%%%%%%%%%%%%%%%%%%%%%%%%%%%%%%%%%%%%%%%%%%%%%%%%%%%%%%

{\bf Proposition 2.} {\it 
The fractional-order variation derivative (\ref{Def-6})
of the functional 
\be
F[u]=\int^{x_2}_{x_1} u(x)u_x(x) dx .
\ee
has the form
\be \label{hx0}
\delta^{\alpha}F[u]=
\int^{x_2}_{x_1} dx \, \left( 
A_1 (\alpha) \, u_x(x) h(x)+ 
A_2 (\alpha) 
u^{2-\alpha}(x) \alpha h^{\alpha-1}(x) h_x (x) \right) ,
\ee
where 
\be
A_1 (\alpha) = 
\left[ \frac{(1-\alpha)\Gamma(2\alpha+1)}{\Gamma(\alpha+1)}+
\frac{ \alpha \Gamma(2\alpha)}{\Gamma(\alpha)} \right] 
\quad
A_2 (\alpha) = \frac{\alpha^2 \Gamma(2\alpha)}{\Gamma(\alpha+1)}
\ee
}

{\bf Proof: }
Using the functional
\[ F[ u^{1-\alpha} (u+\epsilon h)^{\alpha}]=\int^{x_2}_{x_1} dx 
=u^{1-\alpha} (u+ \epsilon h)^{\alpha} 
[u^{1-\alpha} (u+ \epsilon h)^{\alpha} ]_x= \]
\[ =\int^{x_2}_{x_1} dx \Bigl(
(1-\alpha)u^{1-2\alpha} h^{2\alpha} (\epsilon +u/h)^{2\alpha} u_x
+\alpha u^{2-2\alpha} h^{2\alpha-1}(\epsilon+u/ h)^{2\alpha-1} u_x+ \]
\be
 +\alpha u^{2-2\alpha} h^{2\alpha-1} (\epsilon +u/h)^{2\alpha-1} \epsilon h_x 
\Bigr) . \ee
Let us consider $\epsilon_0=\epsilon_0(x)=-u(x)/h(x)$. 
We can use the following relations
\be \label{Frac1}
D^{\alpha}_{\epsilon} (\epsilon-\epsilon_0)^{\beta}=
\frac{\Gamma(\beta+1)}{\Gamma(\beta+1-\alpha)}
(\epsilon-\epsilon_0)^{\beta-\alpha} ,
\ee
and
\be \label{Frac2}
D^{\alpha}_{\epsilon} 
(\epsilon-\epsilon_0)^{\beta} (\epsilon-0)^{\gamma}=
\frac{\Gamma(\beta+1)}{\Gamma(\beta+1-\alpha)}
F_{2;1}\left(-\gamma, \beta+1;\beta+1-\alpha ; 
-\frac{\epsilon -\epsilon_0}{\epsilon_0} \right)
(\epsilon_0)^{\gamma} (\epsilon -\epsilon_0)^{\beta-\alpha} ,
\ee
Note that Eq. (\ref{Frac1}) is satisfyed for $\beta>-1$,
and Eq. (\ref{Frac2}) is satisfyed for $\beta>-1$, 
$\epsilon> \epsilon_0>0$ .
Here $F_{2;1}(a,b,c,z)$ is the Gauss hypergeometric function \cite{SKM}: 
\[ F_{2;1}(a,b,c,z)=\sum^{\infty}_{k=0}
\frac{(a)_k (b_k)}{(c)_k} \frac{z^k}{k!} , \]
and
\[ (z)_k=z(z+1)...(z+n-1)=\frac{\Gamma(z+n)}{\Gamma(z)} . \]
For the function $F_{2;1}(a,b,c,z)$ exists the Euler representation
\be
F_{2;1}(a,b,c,z)=\frac{\Gamma(c)}{\Gamma(b)\Gamma(c-b)} 
\int^1_0 t^{b-1} (1-t)^{c-b-1} (1-zt)^{-a} dt ,
\ee
where $Re(c)> Re (b) >0$, $|arg(1-z)|<\pi$ .
Using
\be
F_{2;1}[a,b;c;1]=\frac{\Gamma(c) \Gamma(c-a-b)}{\Gamma(c-a) \Gamma(c-b)}
\ee
where $Re(c-a-b)>0$ .
For the case $\alpha>1$, we can use
\be
\frac{d^{k}}{dz^k} F_{2;1}[a,b;c;z]=\frac{(a)_k (b_k)}{(c)_k}
F_{2;1}[a+k,b+k;c+k;z] ,
\ee
where $k=[\alpha]$ is a whole part of $\alpha$.
As the result, we get
\[ \delta^{\alpha}F[u]=
\int^{x_2}_{x_1} dx 
\left[ \frac{(1-\alpha)\Gamma(2\alpha+1)}{\Gamma(\alpha+1)}+
\frac{ \alpha \Gamma(2\alpha)}{\Gamma(\alpha)} \right] u_x(x) h(x)+ \]
\be \label{hx}
+\int^{x_2}_{x_1} dx 
\frac{\alpha\Gamma(2\alpha)}{\Gamma(\alpha+1)}
u^{2-\alpha} \alpha h^{\alpha-1} h_x .
\ee
%%% WE use $\Gamma(1-\alpha)/\Gamma(-\alpha)= -\alpha$. Minus?
This end of the proof. \\

For $\alpha \rightarrow 1$, we have the usual relation
\[ \delta F[u]= \int^{x_2}_{x_1} dx [u_xh+uh_x]. \]

Note that we can use
\[ \alpha h^{\alpha-1} h_x= (h^{\alpha})_x, \]
and then 
\be
\int^{x_2}_{x_1} dx 
\frac{\Gamma(2\alpha) \Gamma(1-\alpha)}{\Gamma(-\alpha) \Gamma(\alpha+1)}
u^{2-\alpha} \alpha h^{\alpha-1} h_x 
=-\int^{x_2}_{x_1} dx 
\frac{\Gamma(2\alpha) \Gamma(1-\alpha)}{\Gamma(-\alpha) \Gamma(\alpha+1)}
(u^{2-\alpha})_x h^{\alpha} .
\ee
This allows us to realize the integration by part in second term
of Eq. (\ref{hx}).

%%%%%%%%%%%%%%%%%%%%%%%%%%%%%%%%%%%%%%%%%%%%%%%%%%%%%%%
\section{Fractional Variation of Fields}

In this section we consider the definition of fractional
variational derivative without using the increment 
\[ \Delta F[u]=F[u+h]-F[u]\]
of the functional $F[u]$, 
and without using the derivative with respect 
to parameter $\epsilon$ as in the Gateaux derivative.
We suppose that functional $F[u]$ has some densities $f(u,u_x,...)$. 

If ${\bf u}=(u_1,...,u_m)(x,t)$ be a smooth vector-function
that is defined in the region $W \subset R^n$, then 
the variation of the functional 
\[ F[{\bf u}]=\int_W f({\bf u},{\bf u}_x) dx \]
can be defined by the relation
\be
\delta F[u]= \int_W \delta f({\bf u},{\bf u}_x) dx =
\int_W \left[ \frac{\partial f}{\partial u^{\mu}} \delta u^{\mu}+ 
\frac{\partial f}{\partial u^{\mu}_x} \delta u^{\mu}_x 
\right]
dx .
\ee

To define the fractional generalization of
variation and fractional exterior variational calculus \cite{AK},
we can use an analogy with definition of fractional 
exterior derivative.
If the partial derivatives in the definition of the exterior derivative 
\[ d=dx_i \partial / \partial x_i \] 
are allowed to assume fractional order, 
a fractional exterior derivative can be defined \cite{FDF} 
by the equation
\be
d^{\alpha}=(dx_i)^{\alpha} D^{\alpha}_{x_i} , 
\ee 
where $D^{\alpha}_{x}$ are the fractional derivative with respect to $x$.
Using this analogy, we can define the fractional variation
by the following way. 
For the point $u$ of functional space, 
we can define the fractional variation
$\delta F[u]$ of the functional
\be
F[u]=\int^{x_2}_{x_1} f(u,u_x) dx , 
\ee
where $u_x=du/dx$, by the equation
\be \label{dFu}
\delta^{\alpha} F[u]= 
\int^{x_2}_{x_1} \delta f(u,u_x) dx =
\int^{x_2}_{x_2} dx \left[ (\delta u)^{\alpha}
D^{\alpha}_u f(u,u_x) +
(\delta u_x)^{\alpha} D^{\alpha}_{u_x} f(u,u_x) \right].
\ee
 
This approach has a difficulty with the realization of 
integration by part in the second term of (\ref{dFu}).
It is easy to see that the variation $(\delta u_x)^{\alpha}$
cannot be represented as some operator acts on the variation
$(\delta u)^{\alpha}$, i.e., we have
\be \label{dif2}
\left( \delta u_x(x) \right)^{\alpha}= 
\left(\frac{d}{dx} \delta u(x) \right)^{\alpha} 
\not= \frac{d}{dx} (\delta u)^{\alpha}.
\ee

%%%%%%%%%%%%%%

In order to solve this difficulty we can use the following.

Let us define the fractional variation of the functional
\be \label{Fu1um}
F[u]=\int f(u_1,u_2,...,u_m) dx,
\ee
by the equation
\be \label{FDdua}
\delta^{\alpha} F[u]= \int 
\left[ \sum^{m}_{k=1} (D^{\alpha}_{u_k} f) (\delta u_k)^{\alpha} \right] dx
\ee
by analogy with
\be \label{FDdu1}
\delta F[u]= \int \sum^{m}_{k=1} (D^1_{u_k} f) \delta u_k dx.
\ee

In this case, we have
\be
\delta^{\alpha} u_l= \int 
\left[ \sum^{m}_{k=1} (D^{\alpha}_{u_k} u_l) (\delta u_k)^{\alpha} \right] .
\ee
It is known, that
\be
\delta u_l= \int 
\left[ \sum^{m}_{k=1} (D^1_{u_k} u_l) \delta u_k \right] dx=
\int 
\left[ \delta_{kl} \delta(y-x) \delta u_k \right] dx .
\ee
Using the Caputo fractional derivative, we have
\be
D^{\alpha}_{u_k(y)} u_l(x)= \frac{u^{1-\alpha}}{\Gamma(2-\alpha)} 
\delta_{kl} \delta (y-x) .
\ee
Here the Caputo fractional derivative leads us to the $\delta_{kl}$, i.e.,
\[ D^{\alpha}_{u_k(y)} u_l(x)=0 \quad k\not=l .\]
As the result, we get
\be
\delta^{\alpha} u_l= \frac{u^{1-\alpha}}{\Gamma(2-\alpha)} 
\delta_{kl} (\delta u_k)^{\alpha} .
\ee
Therefore, we have
\be \label{C1}
(\delta u_k)^{\alpha}=\Gamma(2-\alpha) u^{\alpha-1}_k
\delta^{\alpha} u_k .
\ee
For the Riemann-Liouville fractional derivatives,
we can derive some analogous relation:
\be
(\delta u_k)^{\alpha}= A_{kl}\left( \Gamma(2-\alpha) u^{\alpha-1}_k;
\Gamma(1-\alpha) u^{-1}_k u^{\alpha}_l \right) \delta^{\alpha} u_k .
\ee

Substituting Eq. (\ref{C1}) in Eq. (\ref{FDdua}), 
we have the following definition. \\

%%%\be \label{FDdau} \delta^{\alpha} F[u]=\Gamma(2-\alpha) \int 
%%%\left[ \sum^{m}_{k=1} (D^{\alpha}_{u_k} f) u^{\alpha-1}_k 
%%%\delta^{\alpha} u_k \right] dx . \ee

{\bf Definition 7.} {\it 
The fractional order variation $\delta^{\alpha} F[u]$ 
of the functional (\ref{Fu1um}) is defined by the equation
\be \label{FDdau}
\delta^{\alpha} F[u] =\Gamma(2-\alpha) \int 
\left[ \sum^{m}_{k=1} (D^{\alpha}_{u_k} f) u^{\alpha-1}_k 
\delta^{\alpha} u_k \right] dx .
\ee
} \\

We can define the fractional variation for derivative 
$du_k(x)/dx$ by the equation
\be
\delta^{\alpha} \frac{d}{dx} u_k(x)=
\frac{d}{dx} \delta^{\alpha} u_k(x) ,
\ee
where we suppose that the fields $u_k$ do not connected
by some constraint. Analogously, we have
\be
\delta^{\alpha} D^{\beta}_x u_k(x)=
D^{\beta}_x \delta^{\alpha} u_k(x) .
\ee

Let us consider $u_1=u$, and $u_2=u_x$
and the functional
\be \label{uu-x}
F[u]=\int f(u_1,u_2) dx =\int f(u,u_x) dx . \ee
As the result, we have
\be \label{FDuu-x}
\delta^{\alpha} F[u]=
%%%\delta^{\alpha} \int f(u,u_x) dx=
\Gamma(2-\alpha) \int 
\left[ (D^{\alpha}_{u_k} f) u^{\alpha-1}_k 
\delta^{\alpha} u -
\frac{d}{dx} \left( D^{\alpha}_{u_{kx}} f) u^{\alpha-1}_{kx} \right) \right]
\delta^{\alpha} u_k dx .
\ee
We can use this relation only if $u_x$ and $u$ can be considered
as independent values, i.e., for the case
\be \label{indep}
D^{\alpha}_{u_x} u=0, \quad D^{\alpha}_{u} u_x=0 .
\ee

Note that Eq. (\ref{indep}) cannot be satisfied in general case. 
Therefore, we must consider the functional
\be \label{Lagr-mult}
F[u]=\int \left[ f(u_1,u_2) +
\lambda \left(u_2-\frac{d}{dx} u_1 \right) \right] dx . 
\ee
instead of the functional (\ref{uu-x}). \\

{\bf Proposition 3.} {\it 
The fractional-order variation equation 
\[ \delta^{\alpha} F[u]=0 \]
of the functional (\ref{Lagr-mult}) gives
\be \label{FELE}
(D^{\alpha}_{u_1} f) u^{2\alpha-2}_1 -
\frac{d}{dx} \left( u^{\alpha-1}_1 u^{\alpha-1}_2 
D^{\alpha}_{u_2} f \right)=0 .
\ee
where $u_2=du_1/dx$.} \\

{\bf Proof:}
The fractional-order variation of 
the functional (\ref{Lagr-mult}) gives
\[
\delta^{\alpha} F[u] =
\Gamma(2-\alpha) \int dx
\left[ (D^{\alpha}_{u_1} f) u^{\alpha-1}_1 -
\lambda \left( D^{\alpha}_{u_1} \frac{d}{dx} u_1 \right) u^{\alpha-1}_1 
\right] \delta^{\alpha} u_1+ 
\]
\[ 
+ \Gamma(2-\alpha) \int dx
\left[ (D^{\alpha}_{u_2} f) u^{\alpha-1}_2 +
\lambda (D^{\alpha}_{u_2} u_2) u^{\alpha-1}_2 \right]
\delta^{\alpha} u_2 + 
\]
\be
+\Gamma(2-\alpha) \int dx \left[ u_2-\frac{d}{dx} u_1 \right]
(D^{\alpha}_{\lambda} \lambda) \lambda^{\alpha-1}_k 
\delta^{\alpha} \lambda . \ee

Using
\[
\lambda \left( D^{\alpha}_{u_1} \frac{d}{dx} u_1 \right) u^{\alpha-1}_1 
=\lambda u^{\alpha-1}_1 \frac{d}{dx} \left( D^{\alpha}_{u_1} u_1 \right) = 
\]
\be 
=\frac{d}{dx} \left( \lambda u^{\alpha-1}_1 D^{\alpha}_{u_1} u_1 \right)-
D^{\alpha}_{u_1} u_1 \frac{d}{dx} \left( \lambda u^{\alpha-1}_1\right) ,
\ee
we get the following fractional variation of the functional 
\[ 
\delta^{\alpha} F[u]=
\Gamma(2-\alpha) \int dx
\left[ (D^{\alpha}_{u_1} f) u^{\alpha-1}_1 +
D^{\alpha}_{u_1} u_1 \frac{d}{dx} \left( \lambda u^{\alpha-1}_1\right)
\right] \delta^{\alpha} u_1+ 
\]
\[ 
+ \Gamma(2-\alpha) \int dx
\left[ (D^{\alpha}_{u_2} f) u^{\alpha-1}_2 +
\lambda (D^{\alpha}_{u_2} u_2) u^{\alpha-1}_2 \right]
\delta^{\alpha} u_2 + 
\]
\be
+\Gamma(2-\alpha) \int dx \left[ u_2-\frac{d}{dx} u_1 \right]
(D^{\alpha}_{\lambda} \lambda) \lambda^{\alpha-1}_k 
\delta^{\alpha} \lambda . \ee
As the result, we get the field equations
\be \label{1}
(D^{\alpha}_{u_1} f) u^{\alpha-1}_1 +
D^{\alpha}_{u_1} u_1 \frac{d}{dx} \left( \lambda u^{\alpha-1}_1\right)=0,
\ee
\be \label{2}
(D^{\alpha}_{u_2} f)+ \lambda (D^{\alpha}_{u_2} u_2) =0,
\ee
\be \label{3}
u_2-\frac{d}{dx} u_1 =0
\ee
From the Eq. (\ref{2}) ,we derive the Lagrange multiplier
\be
\lambda =-\frac{D^{\alpha}_{u_2} f}{D^{\alpha}_{u_2} u_2} .
\ee
Substituting this equation in Eq. (\ref{1}), we have
\be
(D^{\alpha}_{u_1} f) u^{\alpha-1}_1 -
D^{\alpha}_{u_1} u_1 \frac{d}{dx} \left( 
\frac{u^{\alpha-1}_1}{D^{\alpha}_{u_2} u_2}
D^{\alpha}_{u_2} f \right)=0 .
\ee
Using $D^{\alpha}_u u=u^{1-\alpha} / \Gamma(2-\alpha)$, 
we get (\ref{FELE}).

Equation (\ref{FELE}) is the fractional Euler-Lagrange equation.

\section{Conclusion}

In the general case, the equation of motion
cannot be derived from the stationary action principle.
The class of equations that can be derived from
stationary action principle by using 
fractional variation is a wider class than the usual class 
equations that can be derived by usual (integer, first order) 
variation. The usual equations of motion can be considered as
special case of equations that can be derived by fractional variation 
such that $\alpha=1$.

It seems that the fractional variations are 
abstract and formal constructions. 
For this reason, we would like to pay attention that suggested 
fractional variations can have the wide application in
study of fractional gradient type equations 
and fractional generalization of Lyapunov direct (second)
method in the theory of stability. 

The possible importance of fractional variations 
are connected with following ideas.
The class of gradient dynamical systems is
a restricted class of all dynamical systems.
However these systems have the important property.
The gradient system can be described by one function - potential,
and the study of the system can be reduced to research of potential.
For example, the way of chemical reactions is defined
from the analysis of potential energy surfaces \cite{LB,K1,K2}.
The fractional gradient systems 
has been suggested in Refs. \cite{JPA2005-2,LMP2005}.
The fractional gradient systems are non-gradient dynamical
systems that can be described by one function - 
some potential. 
For example, the Lorenz equations and R\"ossler equations
are fractional gradient systems \cite{JPA2005-2,LMP2005}.
Therefore the study the some non-gradient 
system can be reduced to research of potential.
For example, the way of some chemical reactions with 
dissipation, dynamical chaos and self-organizing 
can be considered by the analysis of some potential energy surfaces.
The suggested fractional variations allows us to define
the fractional generalization of gradient type equations 
that can have the wide application for 
the description dissipative structures \cite{Prig,Zaslavsky7}.
The suggested approach can also be generalized for
lattice systems by using the
lattice fractional calculus \cite{JPA2014}.

%%%%%%%%%%%%%%%%%%%%%%%%%%%%%%
\newpage
\section*{Appendix}

There are the following rules for variational defivatives.

\begin{itemize}

\item The variation derivative of the field $u(x)$ is defined
by the equation
\be \label{uxuy} \frac{\delta u(x)}{\delta u(y)} =\delta(y-x) , \ee
where we use
\be
\delta u(x)=\int \delta(y-x) \delta u(y) .
\ee
The variational derivatives of linear functional 
\be F[u]= \int g(x) u(x) dx \ee 
can be calculated by the simple formula
\[ \frac{\delta F[u]}{\delta u(y)}= 
\int g(x)\frac{\delta u(x)}{\delta u(y)} dx=
\int g(x) \delta (y-x) dx= h(y) . \]

\item If the $u$ function in the functional is affected by 
differential operators, then, 
in order to make use of the rule (\ref{uxuy}), 
one should at first "throw them over" to the left, 
fulfilling integration by parts. For example,
\[ \frac{\delta}{\delta u(y)} \int g(x) [\nabla u(x)] dx =
- \frac{\delta}{\delta u(y)} \int [\nabla g(x)] u(x) dx = -\nabla g(y). \]
We assumed here that on the boundary of integration domain the product
$u(x) h(x)$ becomes zero.

\item 
The variational derivative of nonlinear functionals is calculated according
to the rule of differentiating a complex function similarly to partial 
derivatives:  
\[ \frac{\delta}{\delta u(y)} \int f(u(x)) dx =
\int \frac{\delta f(u)}{\delta u(x)} 
\frac{\delta u(x)}{\delta u(y)} dx = 
\int \frac{\delta f(u)}{\delta u(x)} \delta(y-x) dx = 
\frac{\delta f(u(y))}{\delta u(y)} . \]
For example,
\[ \frac{\delta}{\delta u(y)} \int [u(x)]^n dx = n [u(y)]^{n-1} . \]

\end{itemize}

%%%%%%%%%%%%%%%%%%%%%%%%%%%%%%%%%%%%%%%%%%%%%%%%%%%%%%%%%%%%

\end{document}